\documentclass[11pt]{amsart}
\usepackage{amssymb,amsmath,amsthm}
\usepackage{graphicx}
\usepackage{epsfig}
\newtheorem{theorem}{Theorem}[section]
\newtheorem{lemma}[theorem]{Lemma}

\newtheorem{corollary}[theorem]{Corollary}

\theoremstyle{remark}
\newtheorem{definition}[theorem]{Definition}
\newtheorem{remark}[theorem]{Remark}

\newcommand\Z{{\mathbb Z}}

\begin{document}

\title{Coloring invariants of spatial graphs}
\date{May 24, 2010}
\author{Maciej Niebrzydowski}
\address[Maciej Niebrzydowski]{Department of Mathematics\\
	 University of Louisiana at Lafayette\\
	 1403 Johnston Street\\
	 217 Maxim D. Doucet Hall\\
	 Lafayette, LA 70504-1010}
\email{mniebrz@gmail.com}


\thispagestyle{empty}

\begin{abstract}
We define the fundamental quandle of a spatial graph and several invariants derived from it. In the category of graph tangles, we define an invariant based on the walks in the graph and cocycles from nonabelian quandle cohomology. 
\end{abstract}

\maketitle

\section{Introduction}
In this paper, we define general invariants for spatial graphs that utilize quandles and quandle colorings. There are many such invariants for knots and links 
(see \cite{CKS,FR,Joy,K,N}, for example), and in recent years some of these ideas have been adapted and generalized to graphs (\cite{FM,I,II,IY,MSW}). We will begin with some introductory definitions and examples.

\begin{definition}
A {\it quandle}, $X$, is a set with a binary operation 
$(a, b) \mapsto a * b$
such that
\begin{enumerate}
\item For any $a \in X$, 
$a* a =a$.

\item For any $a,b \in X$, there is a unique $c \in X$ such that 
$a= c*b$. 

\item For any $a,b,c \in X$,
$ (a*b)*c=(a*c)*(b*c)$ (right distributivity). 
\end{enumerate}
\end{definition}

The second condition can be replaced by the following requirement:
the operation $*_b\colon Q\to Q$, defined by $*_b(x)=x*b$, is a bijection. The inverse map to $*_b$ is denoted by
$\overline{*}_b$. It is a standard convention that the operations in a quandle are performed from left to right if the brackets are not used. 

The following examples of quandles often appear in applications:

\begin{itemize}
\item[-]
Any group $H$ with conjugation 
as a quandle operation:\\ $a*b=b^{-1} a b$. Such quandle is called the {\it conjugation quandle} of the group $H$, and denoted by $Conj(H)$.
\item[-]
We can obtain a different quandle from a given group $H$ if we take $a*b=b a^{-1} b$ as a quandle operation. It is called the {\it core quandle} of the group $H$. 
\item[-]
Let $n$ be a positive integer.
For elements  $i, j \in \{ 0, 1, \ldots , n-1 \}$, define
$i\ast j \equiv 2j-i \pmod{n}$.
Then $\ast$ defines a quandle
structure  called the {\it dihedral quandle} and denoted by $R_n$.
It can be identified with  the
set of reflections of a regular $n$-gon
with conjugation
as a quandle operation.
\item[-]
Any $\mathbb{Z}[t, t^{-1}]$-module $M$
is a quandle with
$a*b=ta+(1-t)b$, for $a,b \in M$, called the {\it  Alexander  quandle}.
Moreover, if $n$ is a positive integer, then
$\mathbb{Z}_n[t, t^{-1}]/(h(t))$
is a quandle for
a Laurent polynomial $h(t)$.
\end{itemize}

\begin{definition}
Let $G$ be a graph, and let $f\colon G\to \mathbf{R}^3$ be its embedding. The image $f(G)$ is called a {\it spatial graph}. Two spatial graphs are considered equivalent if one is ambient isotopic to the other. A {\it diagram} of the spatial graph is obtained by projecting $f(G)$ onto a plane such that all self-intersection points are double points that do not include images of vertices. Crossing information is indicated by removing a small part of the arc that was closer to the plane before the projection. During this process the images of edges are split into parts; these parts are called the {\it arcs of the diagram}. 
\end{definition}

\begin{figure}
\begin{center}
\includegraphics[height=12 cm]{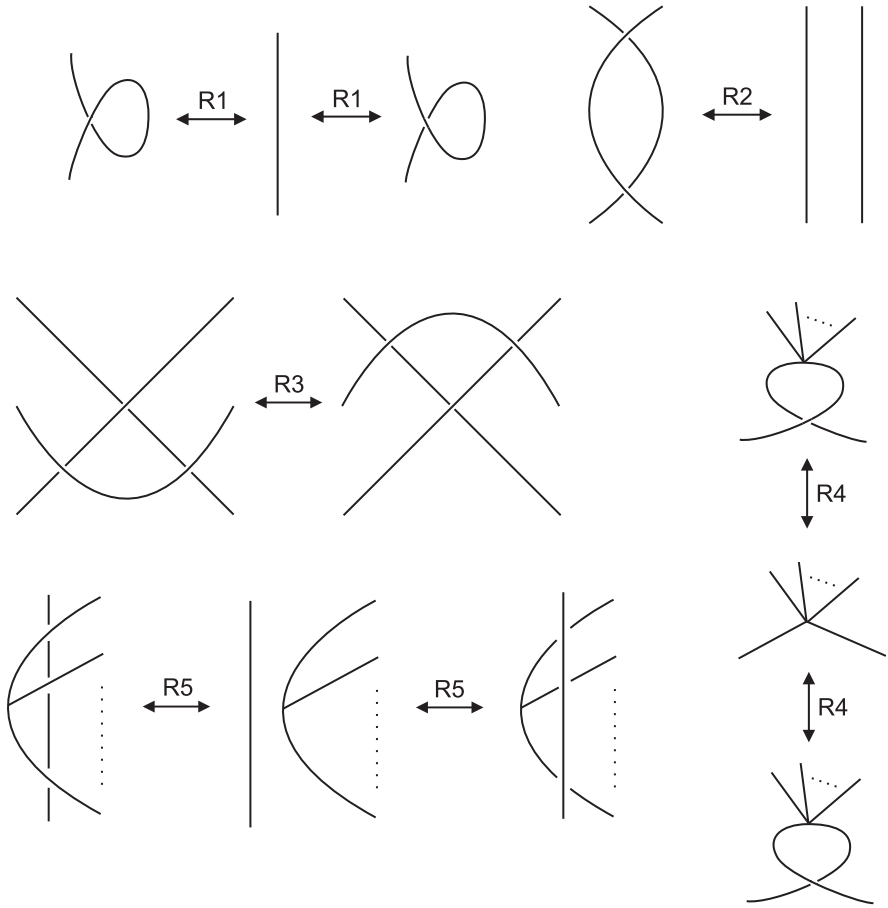}
\caption{Graphical Reidemeister moves.\label{gramoves}}
\end{center}
\end{figure}
 
Ambient isotopy between spatial graphs is described combinatorially by the following theorem \cite{Kauff}.

\begin{theorem}
Two diagrams represent ambiently isotopic spatial graphs if and only if one can be obtained from the other by a finite sequence of graphical Reidemeister moves shown in Figure \ref{gramoves}.
\end{theorem}

\begin{definition}
The spatial graph is {\it oriented} if each of its edges is given a direction. For the oriented spatial graphs we have an analogous theorem to the one above, 
except that we have to use graphical Reidemeister moves with all possible choices of orientation.
\end{definition}

Spatial graph theory is basically a generalization of knot theory (one can add vertices of degree two to the components of a link to obtain a proper spatial graph).
In classical knot theory there is the following notion of a quandle coloring of a link diagram.

\begin{figure}
\begin{center}
\includegraphics[height=3 cm]{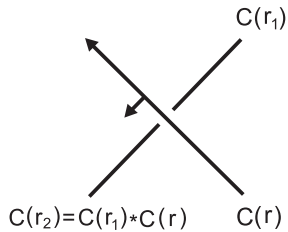}
\caption{The quandle coloring rule for knot and link diagrams.}\label{crossingrule}
\end{center}
\end{figure}

\begin{definition}\label{colorings}
Let $X$ be a fixed quandle, and let $R$ be the set of arcs of the link diagram. The normals to the arcs are given in such a way that the pair (tangent vector, normal vector) matches the usual orientation of the plane. A {\it quandle coloring} is a map $C\colon R\to X$ such that at every crossing, the relation depicted 
in Fig.\ref{crossingrule} is satisfied. More precisely, if $r$ is the over-arc at a crossing, and $r_1,\ r_2$ are under-arcs such that the normal of the over-arc points from $r_1$ to $r_2$, then it is required that 
$C(r_2)=C(r_1)*C(r)$.
\end{definition}

The number of quandle colorings of a link diagram is a link invariant, and, more generally, coloring is the first step in defining many more powerful invariants, including the ones derived from quandle homology and cohomology theories (some of them are described in \cite{CKS}). Several attempts have been made to modify this notion so that it is in harmony with the strictly graphical moves R4 and R5 (see \cite{CKS,FR,Joy,N}, for example). In the next section we will give a general definition of a quandle coloring of a spatial graph that uses the notion of the fundamental quandle.

\section{Homomorphisms from the fundamental quandle}

In this section it will be convenient to use the exponential notation introduced in \cite{FR}. In this notation we write $*_b(a)=a*b$ as $a^b$, and $\overline{*}_b(a)$ as $a^{\overline{b}}$. Now we do not need brackets any more: $a^{bc}$ means $(a^b)^c$ and $a^{b^c}$ means $a^{(b^c)}$. Using simple substitutions the third quandle axiom can be rewritten as $a^{b^c}=a^{\overline{c}bc}$. It emphasizes the connection between the quandle operation and group conjugation, and makes it clear that all expressions can be resolved to one of the form $x^w$, where $x$ is an element of the quandle $X$, and $w$ is a word from the free group on $X$. Elements written in the exponent line are referred to as operators.

\begin{figure}
\begin{center}
\includegraphics[height=3 cm]{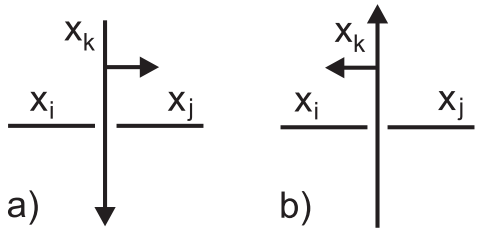}
\caption{Quandle elements assigned to arcs.}\label{crossings}
\end{center}
\end{figure}

\begin{definition}\label{funknot}
The {\it fundamental quandle of an oriented knot} was defined in \cite{Joy}. It is a classifying invariant of classical knots up to the operation that simultaneously changes the orientation of the knot diagram and replaces it by its mirror image. Generators of the fundamental quandle correspond to the arcs of the knot diagram, and relations correspond to crossings. They are of the form: $x_i\ast x_k=x_j$ in the case of a crossing like the one in Fig.\ref{crossings}(a), or $x_i\,\bar{\ast}\,x_k=x_j$, for the crossing like in Fig.\ref{crossings}(b).
\end{definition}

\begin{definition}
For every quandle $X$ there is a corresponding {\it associated group}, $As(X)$, obtained by interpreting the quandle operation as conjugation. More specifically,
$As(X)=F(X)/H$, where $F(X)$ denotes the free group on $X$, and $H$ is a normal subgroup generated by elements $x^yy^{-1}x^{-1}y$, where $x$, $y\in X$.
A group presentation for $As(X)$ can be easily obtained from the quandle presentation of $X$ by taking the same set of generators, and reading the quandle operation in the relations as conjugation (see \cite{FR} for details). For example, $a^b=c^{\overline{d}}$ gives a group relation $b^{-1}ab=dcd^{-1}$.  
It is well known that the associated group of the fundamental quandle of a knot $K$ is isomorphic to the fundamental group of the complement of $K$, and there is a natural action of this group on the fundamental quandle \cite{Joy}. 
\end{definition}

We now proceed to defining the fundamental quandle of an oriented spatial graph.

\begin{figure}
\begin{center}
\includegraphics[height=3.5 cm]{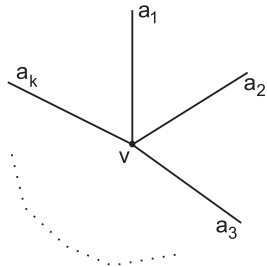}
\caption{Generators assigned to arcs surrounding a vertex.}\label{vertex}
\end{center}
\end{figure}

\begin{definition}\label{main}
Let $D$ be a diagram of an oriented spatial graph $G$, let $V$ be the set of vertices of $D$, and let $R$ denote the set of arcs of $D$. We define the {\it fundamental quandle of $G$} as the quandle with the following presentation:
\begin{itemize}
\item[(1)] To each arc we assign one generator.
\item[(2)] The relations coming from crossings are the same as in the definition of the fundamental quandle of the knot (Definition \ref{funknot}).
\item[(3)] For each vertex $v\in V$, we consider the set of arcs incident to $v$, ordered clockwise. Let $a_1,a_2,\ldots,a_k$ be the set of generators assigned to these arcs (see Fig.\ref{vertex}). Let $\epsilon_i=1$, if the orientation of the arc $a_i$ is pointing towards $v$, and $\epsilon_i=-1$ otherwise. Also, let $s^{t^{\epsilon}}$ denote $s^t$ if $\epsilon=1$, and $s^{\overline{t}}$ if $\epsilon=-1$.
For each generator $x$ of the quandle, we add the relation $x^{a_1^{\epsilon_1}a_2^{\epsilon_2}\ldots a_k^{\epsilon_k}}=x$ to the presentation.   
\end{itemize}
We will use the notation $Q(G)$ for the fundamental quandle of the spatial graph $G$.
\end{definition}

\begin{remark}
By considering all quandle generators in the last step of the above definition, we ensure that $q^{a_1^{\epsilon_1}a_2^{\epsilon_2}\ldots a_k^{\epsilon_k}}=q$ for any quandle element $q$. It follows that the operators $a_1^{\epsilon_1}, a_2^{\epsilon_2}, \ldots, a_k^{\epsilon_k}$ in the relations can be permuted cyclically, and it will not change the quandle. The proof that the isomorphism class of the above quandle is not changed by the moves R1, R2 and R3 is straightforward, the invariance under moves R4 and R5 follows from the fact that the sequences of operators introduced in step (3) correspond to relations in the fundamental group of the complement of the spatial graph (see \cite{Fox} for remarks about such fundamental group), and from the observation that, as operators, elements $s^t$ and $\overline{t}st$ are equal (also, $s^{\overline{t}}$ operates in the same way as $ts\overline{t}$). 
\end{remark}

\begin{remark}
The quandle we defined seems to be more subtle than the quandle defined in \cite{FM}, where the generators assigned to arcs incident to the same vertex were only distinguished by orientation. In fact, the quandle defined there (in the case of classical spatial graphs) is a quotient of the well-involuted double cover of the fundamental quandle $Q(G)$ (see \cite{K} for the definition of the well-involuted double cover of a given quandle).
\end{remark} 

\begin{remark} 
There is an epimorphism of the associated group of $Q(G)$ onto the fundamental group of the complement of $G$. Let us illustrate it with two examples.
\end{remark}
 
\begin{figure}
\begin{center}
\includegraphics[height=3 cm]{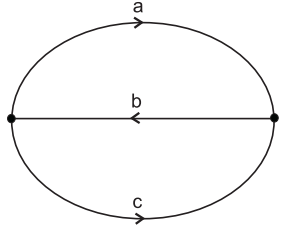}
\caption{An example of a graph with assigned generators.}\label{graphex}
\end{center}
\end{figure}

\noindent\textbf{Example.}
Consider the graph $G$ from Figure \ref{graphex}. The fundamental group of its complement is isomorphic to $\Z * \Z$, as is apparent from its presentation:
$$(a,b,c|a^{-1}bc^{-1}=1, acb^{-1}=1)=(a,b,c|acb^{-1}=1)=(a,b|\,).$$
The presentation of the fundamental quandle $Q(G)$ is as follows:
$$(a,b,c|a^{\overline{a}b\overline{c}}=a,a^{ac\overline{b}}=a,b^{\overline{a}b\overline{c}}=b,b^{ac\overline{b}}=b,c^{\overline{a}b\overline{c}}=c,c^{ac\overline{b}}=c)=$$
$$(a,b,c|a^{b\overline{c}}=a,b^{ac}=b,c^{\overline{b}a}=c).$$
Finally, $As(Q(G))$ is presented as:
$$(a,b,c|cb^{-1}abc^{-1}=a,c^{-1}a^{-1}bac=b,a^{-1}bcb^{-1}a=c).$$ Clearly, the images of relators from $As(Q(G))$ are equal to identity in $\pi_1(\mathbf{R}^3\setminus G)$ when generators are sent to generators. Thus, there is an epimorphism from $As(Q(G))$ onto $\pi_1(\mathbf{R}^3\setminus G)$. 

The fact that this epimorphism is in general not an isomorphism can be clearly seen from the following simple example.

\begin{figure}
\begin{center}
\includegraphics[height=3 cm]{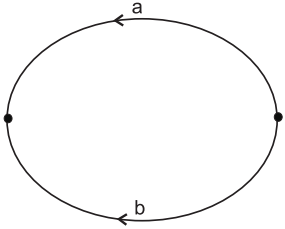}
\caption{A graph with $As(Q(G))$ isomorphic to $\Z\oplus\Z$.}\label{graphex2}
\end{center}
\end{figure}

\noindent\textbf{Example.}
The fundamental group of the complement of the graph from Figure \ref{graphex2} is infinite cyclic:
$$(a,b|ab=1,a^{-1}b^{-1}=1)=(a|\,).$$
Its fundamental quandle has presentation:
$$(a,b|a^{ab}=a, b^{ab}=b, a^{\overline{a}\overline{b}}=a, b^{\overline{a}\overline{b}}=b)=(a,b|a^b=a,b^a=b).$$
That is, $Q(G)$ is a trivial quandle with two elements. Its associated group $As(Q(G))$ is isomorphic to $\Z\oplus\Z$, as can be seen from the presentation:
$$(a,b|b^{-1}ab=a, a^{-1}ba=b)=(a,b|a^{-1}b^{-1}ab=1).$$

\begin{remark}
Let $G$ be a connected plane graph with E edges and V vertices. Then:
\begin{enumerate} 
\item[(i)] The abelianization of $As(Q(G))$ is equal to $\Z^E$.
\item[(ii)] $\pi_1(\mathbf{R}^3\setminus G)$ is a free group of $c(G)= E-V+1$ free generators. 
Thus, its abelianization is equal to $\Z^{E-V+1}$.
\item[(iii)] Equality of abelianizations holds only for one vertex graph (graph of loops).
\end{enumerate}
\end{remark} 

Now we give a general definition of quandle colorings of a spatial graph diagram.

\begin{definition}
Let $X$ be a quandle and let $D$ be a diagram of a spatial graph $G$. We say that an assignment of elements of $X$ to arcs of $D$ is a {\it quandle coloring of a spatial graph diagram $D$} if it is induced from a quandle homomorphism from $Q(G)$ to $X$. More precisely, we consider $Q(G)$ with presentation as in the Definition \ref{main}. If $f$ is a function sending generators of $Q(G)$ to elements of $X$, and the images of relations from $Q(G)$ are satisfied in $X$, then $f$ is a quandle homomorphism, and for each generator $x_i$ of $Q(G)$, we take $f(x_i)$ as the color of an arc to which $x_i$ was assigned.
\end{definition} 

\begin{corollary}
For a fixed quandle $X$, the number of quandle colorings of a diagram $D$ is an invariant of a spatial graph $G$.
\end{corollary}

\begin{definition}
A trivial quandle coloring is a coloring that assigns the same quandle element to every arc of the diagram. Thus, every diagram has at least as many colorings as the size of the quandle that is used.
\end{definition}

\begin{remark}
The invariants defined in \cite{IY} can be understood as the number of nontrivial colorings with the dihedral quandle $R_3$, and the total number of colorings with the dihedral quandle $R_n$.
\end{remark}

\begin{definition}
An $k$-quandle is a quandle $X$ such that $a^{bb\ldots b}=a$ for all $a$ and $b\in X$, where the number of operators $b$ in this condition is $k$.
\end{definition}

\begin{figure}
\begin{center}
\includegraphics[height=6 cm]{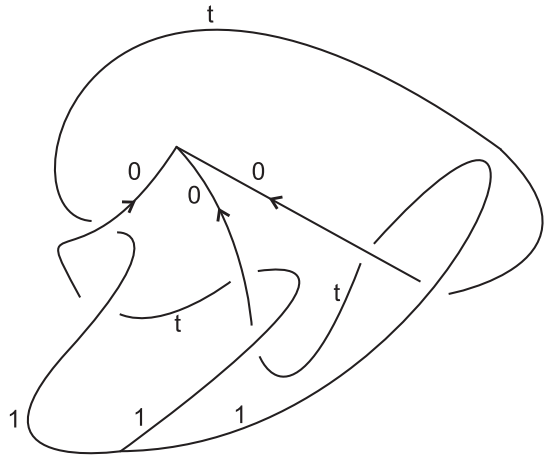}
\caption{Nontrivially knotted theta curve.}\label{theta}
\end{center}
\end{figure}

\begin{remark}
Suppose that $G$ is an oriented spatial graph satisfying the condition that at each of its vertices, the number of incoming edges minus the number of outgoing edges is equal to 0 modulo k. Then, for a given $k$-quandle $X$, it is convenient to consider the kind of colorings in which all arcs incident to the same vertex receive the same color. The number of such special colorings is also a spatial graph invariant. As an example, we show such a coloring for the theta curve depicted in  Figure \ref{theta}, using the Alexander quandle $S_4$, which, as a set, is equal to $\Z_2[t,t^{-1}]/(t^2+t+1)$. The quandle $S_4$ is a $3$-quandle and has four elements $\{0,1,t,t+1\}$. The existence of this nontrivial coloring proves that this theta curve is nontrivially knotted, that is, it is not ambient isotopic to the graph from Figure \ref{graphex} (with any choice of its orientation). 
\end{remark} 

\section{The nonabelian cocycle invariant for graph tangles}


\begin{definition}
By a {\it graph tangle diagram} $\widehat{D}$ we mean a part of a spatial graph diagram $D$ that is contained in a square box, such that the boundary points are not vertices of the graph. A graph tangle $\widehat{G}$ is a class of graph tangle diagrams related by graphical Reidemeister moves R1-R5 and planar isotopy, performed inside the box, with boundary points fixed.
\end{definition}

\begin{definition}
The {\it fundamental quandle of a graph tangle} $\widehat{G}$ is defined in the same way as for a closed spatial graph, by assigning generators to arcs, and relations to crossings and vertices. We also define a {\it quandle coloring of $\widehat{G}$} as an assignment of elements of some fixed quandle $X$ to its arcs that is determined by a quandle homomorphism $Q(\widehat{G})\to X$.
\end{definition}

\begin{definition}
A {\it walk} in a graph is a sequence of vertices and edges $$v_0,e_1,v_1,e_2,\ldots,e_s,v_s,$$ in which each edge $e_i$ joins the vertices $v_{i-1}$ and $v_i$. In the context of graph tangles, we will be interested in {\it extended walks}, the walks that start from a chosen boundary point, and end in the fixed target boundary point. The length of such an extended walk is the number of the whole edges included in it. For oriented graph tangles, it will be useful to call an edge in a given extended walk positive if its original orientation agrees with the direction of the walk, and negative otherwise.
\end{definition}

Now we will define invariants using extended walks, special kind of quandle colorings, and nonabelian cocycles. First, we recall the definition of the second nonabelian quandle cohomology that was defined in \cite{A-G}, and used in \cite{CEGS} to define link invariants.

\begin{definition}
Let $X$ be a quandle and $H$ a group (not necessarily abelian). A function $\phi\colon X\times X\to H$ is a {\it quandle 2-cocycle} if $$\phi(x_1,x_2)\phi(x_1*x_2,x_3)=\phi(x_1,x_3)\phi(x_1*x_3,x_2*x_3)$$
and $$\phi(x,x)=1,$$ for any $x_1$, $x_2$, $x_3$, $x\in X$. The set of quandle 2-cocycles is denoted by $Z_Q^2(X;H)$.
A cocycle $\phi\in Z_Q^2(X;H)$ is called non-abelian when $H$ is not an abelian group. If $H$ is abelian, the above definition gives the usual abelian cocycles, that have been studied in many papers (see \cite{CKS} for literature).
\end{definition}

\begin{definition}
Two cocycles $\phi_1$ and $\phi_2$ are {\it cohomologous} if there is a function $\beta\colon X\to H$ such that
$$\phi_2=\beta(x_1)^{-1}\phi_1(x_1,x_2)\beta(x_1*x_2),$$
for any $x_1$, $x_2\in X$. An equivalence class is called a {\it cohomology class}, and the set of cohomology classes is denoted by $H_Q^2(X;H)$.
\end{definition}

\begin{definition}
Let $\phi\in Z_Q^2(X;H)$, and let $\kappa$ be a crossing with $y_{\kappa}$ being the element of the quandle $X$ assigned to the over-arc, and 
$x_{\kappa}\in X$ being the color of the under-arc from which the normal of the over-arc points. The {\it Boltzmann weight} at $\kappa$ is $\mathcal{B}(\kappa,\mathcal{C})=\phi(x_{\kappa},y_{\kappa})^{\epsilon(\kappa)}$, where $\epsilon(\kappa)$ is $+1$ for the positive crossing $\kappa$, and $-1$ for the negative crossing.
\end{definition} 

In order to define our invariant, we need to impose some conditions. First, we require the cocycle $\phi\in Z_Q^2(X;H)$ to satisfy the equality $$\phi(y*z,z)=\phi(y,z),$$ for any quandle elements $y$, $z\in X$. Notice that when $X$ is finite (and therefore a $k$-quandle, for some $k$), the above condition implies:
$$\phi(y\,\overline{*}\,z,z)=\phi(y(*z)^{k-1},z)=\ldots=\phi(y*z,z)=\phi(y,z).$$  
Secondly, let $s$ be a number such that at each of the vertices of the graph tangle, the number of incoming edges minus the number of outgoing edges is equal to 0 modulo s. Then, we require that $h^s=1$, for all $h\in H$.

Assuming above conditions, let $\mathcal{C}$ be a coloring of an oriented graph tangle diagram $\widehat{D}$ with elements of a finite quandle $X$ such that all arcs incident to the same vertex receive the same color. Let $w$ be an extended walk from a chosen boundary point $b_i$ to another boundary point $b_j$. When going along $w$, we write the Boltzmann weight each time the crossing is encountered in which the followed path is the under-arc, to form the product:
$$\Psi_{(w,\phi,\mathcal{C})}=\phi(x_{\kappa_1},y_{\kappa_1})^{\epsilon(\kappa_1)}\phi(x_{\kappa_2},y_{\kappa_2})^{\epsilon(\kappa_2)}\ldots \phi(x_{\kappa_n},y_{\kappa_n})^{\epsilon(\kappa_n)}.$$

\begin{figure}
\begin{center}
\includegraphics[height=5.5 cm]{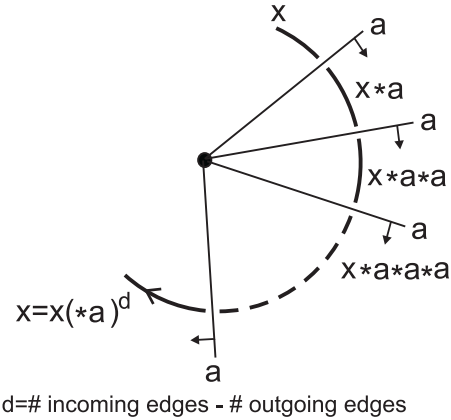}
\caption{Invariance under the move R5.}\label{casest}
\end{center}
\end{figure}

\begin{lemma}
$\Psi_{(w,\phi,\mathcal{C})}$ is not changed by the graphical Reidemeister moves R1-R5.
\end{lemma}
\begin{proof} The invariance under the moves R1 and R4 follows from the equality $\phi(x,x)=1$. The invariance under the second Reidemeister move follows from the fact that the two crossings involved contribute terms differing only by a sign.
The third Reidemeister move does not change $\Psi_{(w,\phi,\mathcal{C})}$ because of the relation $\phi(x_1,x_2)\phi(x_1*x_2,x_3)=\phi(x_1,x_3)\phi(x_1*x_3,x_2*x_3)$. Finally, Figure \ref{casest} depicts one of the possible situations involved in the move R5, where a part of the walk goes below the edges incident to a vertex. Then, the contribution to the product is:
$$\phi(x,a)\phi(x*a,a)\phi(x*a*a,a)\ldots \phi(x(*a)^{d-1},a)=\phi(x,a)^d=1.$$ Here, we used the assumptions listed before.
\end{proof}

Let $\mathcal{C}_q$ be the family of all colorings of $\widehat{D}$ such that all arcs incident to the same vertex receive the same color, and the color of the arc originating from the point $b_i$ is equal to a fixed $q\in X$. Let 
$\mathcal{W}=\mathcal{W}(b_i,b_j,n)$ be the set of all extended walks from $b_i$ to $b_j$ that have length not greater than $n$.
Given $\phi$ and $H$ as before, we can define two multisets:
$$\Psi_q(\widehat{D})=\{\Psi_{(w,\phi,\mathcal{C})}\}_{w\in\mathcal{W}, \mathcal{C}\in\mathcal{C}_q}$$
and
$$\Psi(\widehat{D})=\{\Psi_q(\widehat{D})\}_{q\in X}.$$

\begin{theorem}
$\Psi_q(\widehat{D})$ and $\Psi(\widehat{D})$ are not changed by the graphical Reidemeister moves R1-R5.
\end{theorem}
\begin{proof}
The proof follows from the fact that each $\Psi_{(w,\phi,\mathcal{C})}$ is not changed by the moves R1-R5. 
\end{proof} 

\begin{remark}
\begin{enumerate}
\item[(1)] In some cases it might be useful to limit the above invariant by assuming that the walks used are of special kind. For example, we could use paths (no vertex is repeated), or trails (no edge is repeated), or choose only walks with a fixed length and/or specified number of positive edges.
\item[(2)] The above invariant can be calculated for any pair of boundary points, and the collection of values for different pairs is another invariant. If all the possible pairs are considered, then the order of boundary points can be disregarded.
\item[(3)] If the graph tangle satisfies the condition that at each vertex the number of incoming edges equals the number of outgoing edges, then the restrictions on the group $H$ and the additional cocycle condition $\phi(y*z,z)=\phi(y,z)$ are not necessary.
\end{enumerate} 
\end{remark} 

\begin{figure}
\begin{center}
\includegraphics[height=4.5 cm]{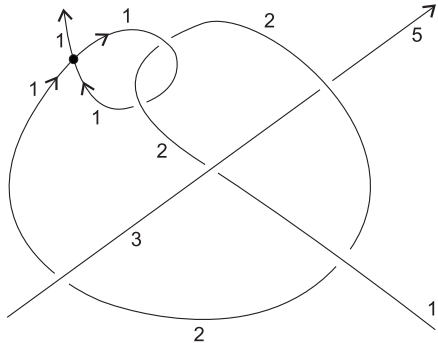}
\caption{A tangle graph with nontrivial cocycle invariant.}\label{przgr}
\end{center}
\end{figure}

\noindent\textbf{Example.}
The following is the multiplication table for the quandle listed on the page 177 of \cite{CKS} that has rich quandle cohomology for $H$ being the symmetric group on three elements.

\begin{center}
\begin{tabular}{c| c ccccc} 
$*$&1 & 2 & 3 & 4 & 5  \\
\hline 
1 & 1 & 1 & 2 & 2 & 2  \\
2 & 2 & 2 & 1 & 1 & 1  \\
3 & 4 & 5 & 3 & 5 & 4  \\
4 & 5 & 3 & 5 & 4 & 3  \\
5 & 3 & 4 & 4 & 3 & 5  \\
\end{tabular}
\end{center}

\noindent We found using GAP \cite{GAP} that there are 8 nontrivial cocycles in different cohomology classes that satisfy the condition 
$\phi(y*z,z)=\phi(y,z)$. Here is one of these cocycles:\\
$\phi(2,1)=(1,2,3)$, $\phi(1,2)=(1,3,2)$, $\phi(1,3)=\phi(2,3)=(2,3)$, $\phi(1,4)=\phi(2,4)=(1,2)$, $\phi(1,5)=\phi(2,5)=(1,3)$, $\phi=()$ for the remaining pairs of quandle elements.

Figure \ref{przgr} shows an example of a graph tangle diagram $\widehat{D}$ that has nontrivial colorings with the above quandle.
In order to distinguish this tangle from its mirror image, let us consider all walks of length zero (recall that only entire graph edges contribute to the length) from the lower right to the upper left boundary point. We will only use the colorings with the color of the lower right boundary point equal to 1. There are five of them. The contribution from the coloring shown in the Figure is:
$$\phi(1,3)\phi(2,1)\phi(2,1)\phi(2,3)=(1,2,3).$$
When all five colorings are taken into account, we obtain the invariant
$$\Psi_q(\widehat{D})=\Psi_1(\widehat{D})=\{(),(1,2,3),(1,2,3),(1,2,3),(1,2,3)\}.$$
The invariant for the mirror image is:
$$\Psi_1(\overline{D})=\{(),(1,2,3),(1,2),(1,3),(2,3)\}.$$

\section{Acknowledgements}
M.~Niebrzydowski is being supported by the Louisiana Board of Regents grant (\# LEQSF(2008-11)-RD-A-30).
The author has had valuable conversations with J. Scott Carter, J\'{o}zef H. Przytycki and Masahico Saito about this paper.

\end{document}